\documentclass{article}

\usepackage{amstext,amsfonts,graphicx,latexsym}

\def\E{\mathop{\hbox{\sf E}}\nolimits}
\def\P{\mathop{\hbox{\sf P}}\nolimits}

\def\phi{\varphi}

\def\B{\Big}

\def\bs{\backslash}

\def\R{{\mathbb R}}
\def\Z{{\mathbb Z}}
\def\C{{\mathbb C}}
\def\H{{\mathbb H}}

\newtheorem{lemma}{Lemma}[section]
\newtheorem{theorem}{Theorem}

\usepackage{graphics}

\def\eqref#1{(\ref{eq.#1})}

\def\putfigure#1#2{
	\begin{figure}[ht]
	\centering
	\includegraphics{#1.eps}
	\caption{#2}
	\label{fig.#1}
	\end{figure}
}

\def\defined#1{{\em #1}}

\def\optional#1{}

\def\proof{\par\noindent{\bf Proof.\ }}

\def\eop{\vskip 3mm }

\newcommand{\leb}{{\mathcal L}}

\begin{document}

\title{Connected allocation to Poisson points in $\R^2$}
\author{Maxim Krikun}
\maketitle
\begin{abstract}
This note answers one question in \cite{hoffman-2006-34},
concerning the connected allocation for the Poisson process in $\R^2$.
The proposed solution makes use of the Riemann map from the plane minus 
the minimal spanning forest of the Poisson point process to the halfplane.
A picture of numerically simulated example is included.
\end{abstract}

\section{The problem}
Let $X\subset \R^d$ be a discrete set. 
We call the elements of $X$ \defined{centers},
the elements of $\R^d$ -- \defined{sites}, 
and we write $\leb$ for the Lebesgue measure in $\R^2$.
%
An \defined{allocation} of $\R^d$ to $X$ 
with \defined{appetite} $\alpha\in[0,\infty]$
is a measurable function $\psi:\R^d \to X \cup\{\infty,\Delta\}$
such that $\leb\psi^{-1}(\xi)\leq\alpha$ for all $\xi\in X$ 
and $\leb\psi^{-1}(\Delta)=0$.
We call $\psi^{-1}(\xi)$ the \defined{territory} of the center $\xi$.
A center $\xi$ is \defined{sated} if $\leb\psi^{-1}(\xi)=\alpha$
and \defined{unsated} otherwise. 
A site $x$ is \defined{claimed} if $\psi(x)\in X$
and \defined{unclaimed} if $\psi(x)=\infty$.
The allocation is \defined{undefined} at $x$ if $\psi(x)=\Delta$.

One particular question we're interested in is the following \cite{hoffman-2006-34}:
{\em
Is there a translation-equivariant allocation 
for the Poisson process of unit intensity
in the critical two-dimensional case ($d=2$ and $\alpha=1$),
in which every territory is connected?
}
The goal of the present note is to describe one such allocation.
Note: for a far more general approach to this question in $d\ge3$
 see \cite{chatterjee-2006}.

\section{Construction of a connected allocation}

Let $X$ be a realization of Poisson process (of intensity $1$) in $\R^2$.

Let the \defined{minimal spanning forest} $T$ be the union of edges
$e=(x,y)$, $x,y\in X$, such that there is no path from $x$ to $y$
in the complete graph on $X$ with all edges strictly shorter than $e$
(see \cite{lyons-2006-34} for other definitions and references on the subject).

For a 2-dimensional Poisson process Alexander \cite{alexander-1995-23}
proved that $T$ is a.s. a tree with one topological end
(a topological end in a graph is a class of equivalence of semi-infinite
 paths modulo finite symmetric difference).
Thus the domain $D = \R^2\backslash T$ is connected, and 
by the Riemann mapping theorem 
can be mapped conforamlly to the upper half-plane $\H = \{ z\in \C: \Im(z)>0 \}$.
Moreover, we can choose such a mapping $f$ which also sends infinity to infinity.
(Note that 
any two such functions differ by a conformal automorphism of $\H$
which invariates infinity, i.e. by an affine mapping 
$z\to az+b$, $a,b\in\R$, $a\neq0$).


The mapping $f$ cannot be extended unambiguously to include $T$.
Instead, consider the inverse mapping $g=f^{-1}$ from the open halfplane
to $D$, then extend $g$ to the real axis by continuity,
obtaining a surjection $\overline{g}: \overline\H\to\R^2$.
Some points of $X$ may have more than one preimage under $\overline{g}$,
more exactly, a point has as many preimages as it's degree in $T$.
For those points we'll have to split their appetite between these preimages
in some way, e.g. proportionally to the angle of the corresponding sector.
Let $(x'_n,\alpha'_n)$ be the resulting points and appetites.
In the following we call the countable set $\{ (x'_n,\alpha'_n), n\in \Z \}$ 
\defined{the image of $X$ under $f$}.

\putfigure{fig2}{The image of $X$ under $f$;
the unit appetite of the center $A$ 
is split into three unequal parts.}


\noindent
Now consider the following \defined{stable allocation procedure}
(see Appendix~II for a formal definition):
\begin{itemize}
\item
Each center starts growing a ball centered in it.
All the balls grow simultaneously, at the same linear speed.
%
\item
Each center claims the sites captured by it's ball, 
unless the site was claimed earlier by some other center.
\item
Once the center becomes sated (the measure of it's territory
reaches it's appetite), the ball stops growing.
\end{itemize}
Note that from our choice of function $f$ it follows that the set
$\{ x'_n, n\in\Z\}$ is locally finite (i.e. has no accumulation point
other than infinity), so the stable allocation procedure is well defined.

Applying this procedure to the image of $X$ under $f$ 
in the half-plane $\H$, with respect to the Eucledian metric in $\H$
and the image $\lambda$ of the Lebesgue measure $\leb$ under $f$
yields certain allocation $\psi_H: \H\to X\cup\{\infty,\Delta\}$.

\begin{lemma}
The allocation $\psi_H$ thus constructed
has connected territories and 
is invariant under affine transformations of $\H$.
\end{lemma}
\proof
The affine invariance follows immediately by construction.

To show that each center has connected territory, consider the following.
Let $A$ be a center, necessarily located at the boundary of $\H$,
and let $D_A = \psi_H^{-1}(A)$ be it's territory.
Consider the ray $AA'$, perpendicular to $\R$.
If some point $x$ on $AA'$ doesn't belong to $D_A$,
that's because the center got sated before reaching $x$,
so no further point on $AA'$ belongs to $D_A$.
Thus the intersection of $D_A$ with $A$ is a segment $AK$.

\putfigure{fig1}{Every point of $D_A$ can be reached from $A$}

Now consider an arc of a circle centered in $A$ and intersecting $AK$ 
in some point $Y$, and follow this arc to the right (or left) from $Y$.
Once we meet a point $Z\notin D_A$, 
there must be a center $B$ such that $BZ\le AZ$,
so no further points on the arc belong to $D_A$.
Thus every point of $D_A$ can be reached from $A$, and $D_A$ is connected.
\eop

Now apply the inverse map $f^{-1}$ to $\psi_H$.
i.e. let $\psi = f^{-1} \psi_H f$.
Clearly, $\psi$ is an allocation of $\R^2$ to $X$ with appetite $\alpha$,
and from the previous lemma it doesn't depend on the choice of $f$.

\begin{lemma}
The allocation $\psi$ is translation-invariant,
i.e. if $\tau: \R^2 \to \R^2$ is a translation,
and $\psi'$ it the allocation or $\R^2$ to $\tau X$ 
constructed using the above procedure, then 
$\psi' = \tau^{-1} \psi \tau$.
\end{lemma}
\proof
Let $X'=\tau X$.
Clearly, the minimal spanning forest is translation-invariant,
i.e. $T' := MSF(X') = \tau \cdot MSF(X)$.
Now $f' = f\,\tau^{-1}$ is a conformal function 
that maps $\R^2\backslash T'$ to $\H$,
and the image of $X'$ under $f'$ coincides with the image of $X$ under $f$,
thus $\psi' = f'^{-1}\psi_H f' = \tau^{-1} f^{-1} \psi_H f \tau = \tau^{-1}\psi\tau$.
\eop

\begin{lemma}
Under the allocation $\psi$ every center is sated a.s.
\end{lemma}

\proof
First, it follows from the Galey-Shapley algorithm (see Appendix~II below), that 
{\em no stable allocation may have both unclaimed sites and unsated centers.}
By ergodicity, the existence of unclaimed sites is a $0/1$ event;
thus we may assume that $\psi_H$, and therefore $\psi$, have no unclaimed sites.

Now Lemma~16 in \cite{hoffman-2006-34} states that for any translation-invariant 
allocation $\psi$ and any $r>0$
\[ \P\{ |\psi(0)| < r \} = \E^* \leb \B( \psi^{-1}(0) \cap B(0,r) \B), \]
where $\E^*$ denotes the expectation with respect to 
the Poisson process conditioned to have a center in $0$.
Taking $r\to\infty$ yields
\[ \P\{ \mbox{$0$ is clamied}\} = \E^* \leb \psi^{-1}(0). \]
But we assumed the probability on the left to be equal to $1$,
and $\leb\psi^{-1}(0) \le 1$ by construction (a center cannot allocate more than
it's appetite), thus the lemma follows.
\eop
\noindent We summarize the above lemmas in the following 
\begin{theorem}
The allocation $\psi$ thus constructed has connected territories, 
is translation-invariant,
and every center is sated a.s.
\end{theorem}

\section{Further questions}
{\em Geometry of territories.} 
Is it true that the every territory $\psi^{-1}(\xi)$ is bounded a.s.?
%
%
The territories of $\psi$ are connected, but, a priori, may be 
not simply connected (more exactly, their closures may not be simply connected).
Can the construction of $\psi$ be modified to assure this?

{\em Point process on $\R$.}
What can be told about the image of $X$ under the Riemann map,
as a point process on $R$?
Simulations suggest that this process should have extremely high
variation of point density; even for small polygons of 30-40 points,
the distances between the images of consecutive points may vary from 
$\approx 1$ to $\approx 10^{-16}$.

{\em Fixing the Riemann map.}
The Riemann map $f$ is unique up to a conformal automorphism of $H$,
$z\to az+b$. 
Is it possible to pick (deterministically) one particular mapping from this class?
If omitting the translation-invariance, one obvious way to do this 
is to pick a mapping $f$ that sends $0$ to $i$.
If doing this in translation-invariant manner, there is obviously no
way to fix the shift $b$ (since this would imply the deterministic
translation-invariant choice of single point from the Poisson process).
But does there exist a way to pick the {\em scale}, i.e. fix the parameter $a$?

\bibliographystyle{hplain} 
\bibliography{alloc}

\newpage
\section*{Appendix I: An obligatory picture}
\begin{figure}[ht]
  \centering
  \includegraphics[width=105mm]{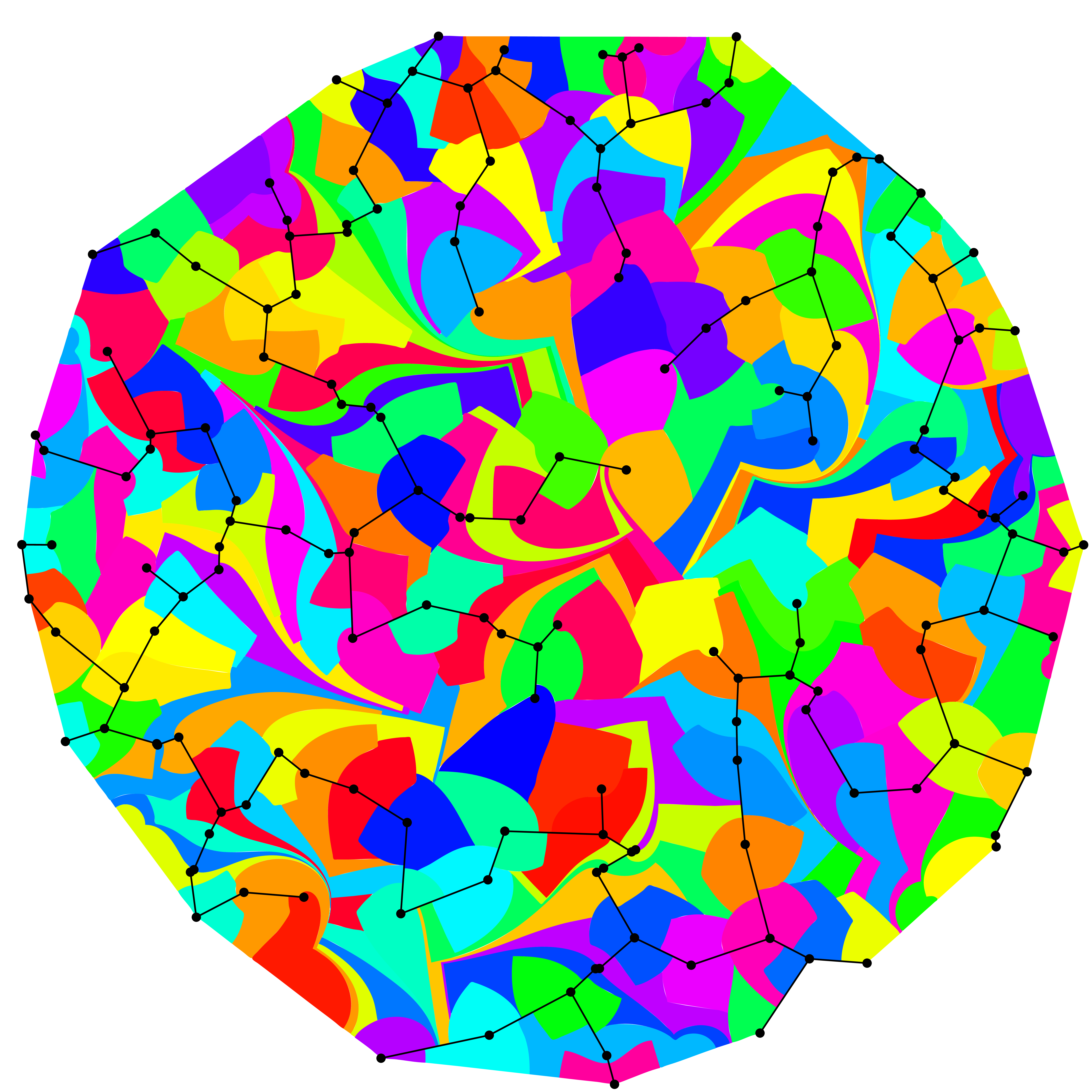}
  \caption{An approximation of allocation for 164 points in the unit circle}
  \label{fig.alloc}
\end{figure}
On the above picture we took 164 uniformly distributed points inside the unit circle,
and mapped each component $K_j$ of their convex hull minus their MST
to a region of $\H$, so that all the internal edges are mapped to $\R$.
Then for each component we allocated it's whole area to the bounding corners,
choosing the appetites proportional to the angular measure.

The mapping was approximated with a version of the {\em geodesic algorithm},
adapted for domains with non-Jordan boundary (or more exactly, for domains 
with boundary consisting of a Jordan curve plus a number of internal trees);
for discussion of numeric algorithms for quasi-conformal mappings 
see e.g.~\cite{marshall-2006} and references therein.

Since the ratio area/\{sum of appetites\} varies in different components $K_j$,
not all domains have exactly the same area;
there are also some artifacts related to the limited numeric accuracy.
Nevertheless, we expect the picture close to the center to be an accurate
approximation of the construction described in this note.

\section*{Appendix II: site-optimal Gale-Shapley algorithm}
{\em 
(The description of the G-S algorithm below is taken from \cite{hoffman-2006-34}
except for the last step, when we have to be more careful with the definition 
of $\psi$ on $W$ in order to assure the connectedness of the territories.)}

Let $W$ be the set of all sites, equidistant from one or more centers.
Since the set of centers is countable, $W$ has Lebesgue measure null.

We construct $\psi$ on $\R^2 \bs W$ by means of a sequence of stages, where 
stage $n$, $n=1,2,3,\ldots$ consists of two steps:
\begin{itemize}
\item[a)] Each site $x\notin W$ \defined{applies} to the closest center,
which has not \defined{rejected} $x$ at any earlier stage.
\item[b)] For each center $\xi$, let $A_n(\xi)$ be the set of sites which
applied to $\xi$ on step a of stage $n$, and define the \defined{rejection radius} as
\[ r_n(\xi) = \inf\{ r: \leb (A_n(\xi) \cap  B(\xi,r)) \ge \alpha \}, \]
where the infimum over the empty set is taken to be $\infty$.
Then $\xi$ \defined{shortlists} all sites in $A_n(\xi) \cap B(\xi, r_n(\xi))$,
and \defined{rejects} all sites in $A_n(\xi)\backslash B(\xi, r_n(\xi))$.
\end{itemize}

Now either $x$ is rejected by every center (in order of increasing distance from~$x$),
or $x$ is shortlisted by $\xi$ for some center $\xi$ at some stage $n$.
In the former case we put $\psi(x)=\infty$ (so $x$ is unclaimed),
in the latter case we put $\psi(x)=\xi$.

Finally, for $x\in W$ put $\psi(x)=\xi$, if $\xi$ is the only center
such that $x \in \partial \psi^{-1}(\xi)$,
and put $\psi(x) = \Delta$ (undefined) otherwise.

\vskip 3cm
\noindent
Maxim Krikun \\
Institut Elie Cartan, \\
Universite Henri Poincare,\\
Nancy, France\\
krikun@iecn.u-nancy.fr \\

\end{document}